\documentclass[1pt,notitlepage,twoside,a4paper]{amsart}

\usepackage{amsmath,amssymb,enumerate}

\usepackage{epsfig,fancyhdr,color}%,showkeys,amsmidx
\usepackage{amssymb}
\usepackage{amsmath,amsthm}    
\usepackage{latexsym}
\usepackage{amscd} 
\usepackage{psfrag}
\usepackage{graphicx} 
\usepackage[latin1]{inputenc}  
\usepackage[all]{xy} 
\usepackage[mathcal]{eucal}

\theoremstyle{definition}

\theoremstyle{remark}

\def\interieur#1{\mathord{\mathop{\kern 0pt #1}\limits^\circ}}

\definecolor{NoteColor}{rgb}{1,0,0}

 \title[A note on Nicolas-Auguste Tissot]{A note on Nicolas-Auguste Tissot: 
At the origin of quasiconformal mappings}

\author{Athanase Papadopoulos}
\address{Athanase Papadopoulos,  Universit{\'e} de Strasbourg and CNRS,
7 rue Ren\'e Descartes,
 67084 Strasbourg Cedex, France}
\email{athanase.papadopoulos@math.unistra.fr}

 \date{\today}

\makeindex
\begin{document}

\begin{abstract}  
Nicolas-Auguste Tissot (1824--1897) was a French mathematician and cartographer. He introduced a tool which became known among geographers under the name \emph{Tissot indicatrix}, and which was widely used during the first half of the twentieth century in cartography. This is a graphical representation of a field of ellipses, indicating at each point of a geographical map the distorsion of this map, both in direction and in magnitude. Each ellipse represented at a given point is the image of an infinitesimal circle in the domain of the map (generally speaking, a sphere representing the surface of the earth) by the projection that realizes the geographical map. 

Tissot studied extensively, from a mathematical viewpoint, the distortion of mappings from the sphere onto the Euclidean plane, and he also developed a theory for the distorsion of mappings between general surfaces. His ideas are close to those that are at the origin of the work on quasiconformal mappings that was developed several decades after him by Gr\"otzsch, Lavrentieff, Ahlfors and Teichm\"uller. 

Gr\"otzsch, in his papers, mentions the work of Tissot, and in some of the drawings he made for his articles, the Tissot indicatrix is represented.
 Teichm\"uller mentions the name Tissot in a historical section in one of his fundamental papers in which he points out that quasiconformal mappings were initially used by geographers. 
 
  The name Tissot is missing from all the known historical reports on quasiconformal mappings.
 In the present article, we report on this work of Tissot, showing that the theory of quasiconformal mappings has a practical origin.
 
 The final version of this article will appear in Vol. VII of the Handbook of Teichm\"uller Theory (European Mathematical Society Publishing House, 2020).
  
\bigskip



\end{abstract}
\larger

\maketitle

  \tableofcontents
 \section{Introduction}
 
Darboux, starts his 1908 ICM talk whose title is \emph{Les origines, les m\'ethodes et les probl\`emes de la g\'eom\'etrie infinit\'esimale} (The origins, methods and problems of infinitesimal geometry)  with the words: ``Like many other branches of human knowledge, infinitesimal geometry was born in the study of practical problems,"  and he goes on explaining how problems that arise in the drawing of geographical maps, that is, the
representation of regions of the surface of the Earth on a Euclidean piece of paper, led to the most important developments in geometry made by Lagrange, Euler, Gauss and others. 
  
\begin{figure}[htbp]
\centering
\includegraphics[width=5.9cm]{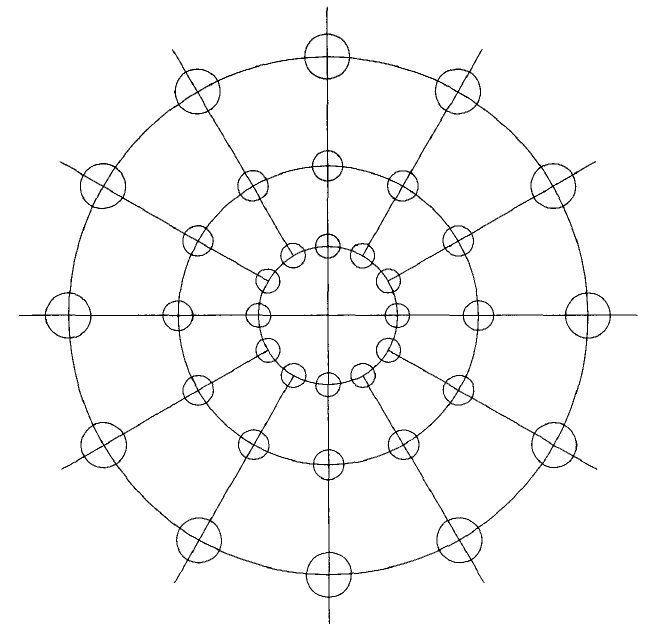}
\includegraphics[width=5.6cm]{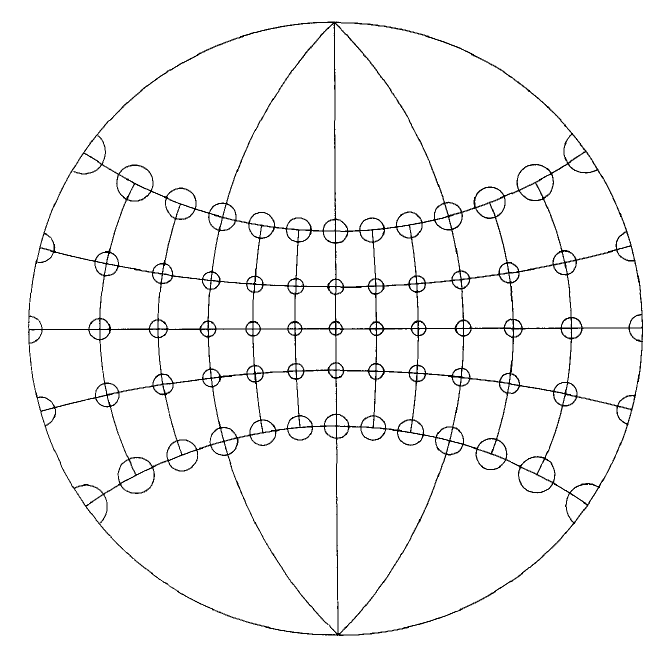}
\includegraphics[width=5.5cm]{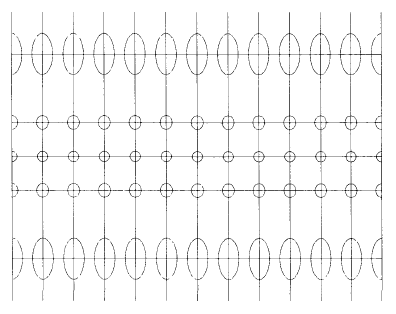}
\includegraphics[width=6.2cm]{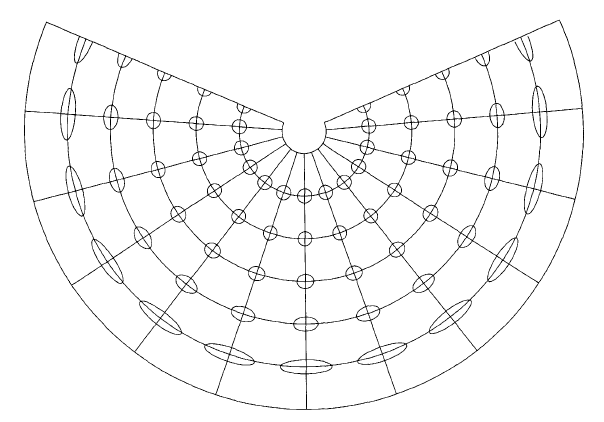}
\caption{\small{Four geographical maps on which  the field of ellipses (Tissot indicatrix) are drawn. The maps are extracted from the book \emph{Album of map projections} \cite{Snyder-Album}. These are called, from left to right, top to bottom, the \emph{stereographic} \cite[p. 180]{Snyder-Album}, \emph{Lagrange} \cite[p. 180]{Snyder-Album},  \emph{central cylindrical} \cite[p. 30]{Snyder-Album} and \emph{equidistant conical} projections \cite[p. 92]{Snyder-Album}.  The first two projections are conformal and not area-preserving. The last two are neither conformal nor area-preserving.}} \label{Snyder}
\end{figure}

The theory of quasiconformal mappings\index{quasiconformal mapping} has its origin in the problems of drawing geographical maps.\index{geographical map}  
Teichm\"uller, in the last part of his paper \emph{ Extremale quasikonforme Abbildungen und quadratische Differentiale} (Extremal quasiconformal mappings and quadratic differentials), published in 1939 \cite{T20}, which is the main paper in which he develops the theory that became known as \emph{Teichm\"uller theory}, makes some comments on this origin, mentioning the work of the French mathematician and geographer Nicolas-Auguste Tissot  (1824--1897).\index{Tissot, Nicolas-Auguste (1824--1897)} Gr\"otzsch,\index{Gr\"otzsch, Herbert (1902--1993)}  in his paper \emph{\"Uber die Verzerrung bei nichtkonformen schlichten Abbildungen mehrfach
zusammenh\"angender schlichter Bereiche} (On the distortion of non-conformal schlicht mappings of multiply-connected schlicht
regions), published in 1930 \cite{Groetzsch1930}, mentions several times the name Tissot, referring to the \emph{Tissot indicatrix}\index{Tissot indicatrix}\index{indicatrix!Tissot} which he represents in the pictures he drew for his article. The directions of the major and minor and minor axes of this ellipse constitute are important element in some of his results. A geographical map\index{geographical map}  is the image of a mapping---henceforth called a projection---from the surface of the Earth, considered as a sphere or spheroid, onto the Euclidean plane.  The Tissot indicatrix is a device introduced by Tissot, who called it the \emph{indicating ellipse} (ellipse indicatrice, which was used by geographers until the middle of the twentieth century. It is a field of ellipses drawn on the  geographical map,\index{geographical map}  each ellipse representing the image by the projection---assumed to be differentiable---of an infinitesimal circle\footnote{The expression ``infinitesimal circle" means here, as is usual in the theory of quasiconformal mappings, a circle on the tangent space at a point. In practice, it is a circle on the surface which has a ``tiny radius." In the art of geographical map drawing, these circles, on the domain surfaces, are all supposed to have the same small size, so that the collection of relative sizes of the image ellipses becomes also a meaningful quantity.} at the corresponding point on the sphere (or spheroid) representing the surface of the Earth. Examples of Tissot indicatrices\index{Tissot indicatrix}\index{indicatrix!Tissot} are given in Figure \ref{Snyder}.

\begin{figure}[htbp]
\centering
\includegraphics[width=6cm]{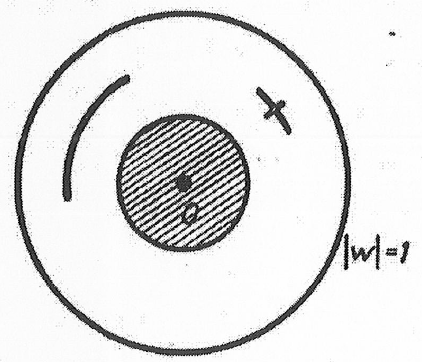}
\includegraphics[width=6.2cm]{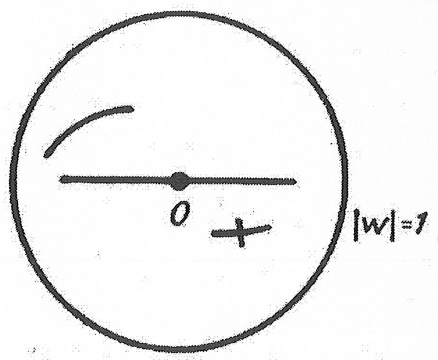}
\caption{\small{Two figures from Gr\"otzsch's paper \emph{\"Uber die Verzerrung bei nichtkonformen schlichten Abbildungen mehrfach
zusammenh\"angender schlichter Bereiche} \cite{Groetzsch1930}.
Gr\"otzsch drew the Tissot indicatrices of his quasiconformal mappings. (In each drawing, the major and minor axes of the ellipses are shown.)}} \label{Groe}
\end{figure}

  In Figure \ref{Groe}, we have reproduced drawings from a paper of Gr\"otzsch\index{Gr\"otzsch, Herbert (1902--1993)}  in which he represents the Tissot indicatrix\index{Tissot indicatrix}\index{indicatrix!Tissot} of the maps he uses.
  
Although the work of Tissot\index{Tissot, Nicolas-Auguste (1824--1897)} is closely related to the theory of quasiconformal mappings, his name is never mentioned in the historical surveys of this subject, and the references by Gr\"otzsch and by Teichm\"uller to his work remained unnoticed. In this note, I will give a few indications on this work.

Before surveying the work of Tissot\index{Tissot, Nicolas-Auguste (1824--1897)} in \S \ref{s:work}, I will give, in \S \ref{s:bio}, a short biographical note on him.

 \section{Biographical note on Tissot}\label{s:bio}

Nicolas-Auguste Tissot was born in 1824,\index{Tissot, Nicolas-Auguste (1824--1897)} in  Nancy, which was to become, 26 years later, the birthplace of Henri Poincaré (whom we shall mention soon). Tissot entered the \'Ecole Polytechnique in 1841. He started by occupying a career in the Army\footnote{The reader should note that the \'Ecole Polytechnique was, and is still, is a military school.}  and defended a doctoral thesis on November 17, 1851; cf. \cite{Tissot-these}. On the cover page of his thesis, he is described as ``Ex-Capitaine du G\'enie."  Tissot became later a professor at the famous Lycée Saint-Louis in Paris, and at the same time examiner at the \'Ecole Polytechnique, in particular for the entrance exam. He eventually\index{Tissot, Nicolas-Auguste (1824--1897)} became an assistant professor (\emph{r\'ep\'etiteur}) in geodesy at the  \'Ecole Polytechnique. 

After having published, in the period 1856--1858, several papers and \emph{Comptes Rendus} notes on cartography,\index{geographical map} in which he analyzed the distortion of some known geographical maps\index{geographical map}  (see \cite{Tissot-CR00}, \cite{Tissot-1858}, \cite{Tissot-CR0}), Tissot started developing his own theory, on which he published  three notes, in the years 1859--1860,  \cite{Tissot-CR1, Tissot-CR2, Tissot-CR3}, and then a series of others in the years 1865--1880  \cite{Tissot-CR4, Tissot1878, Tissot1878a, Tissot1878b, Tissot1879, Tissot1879a, Tissot1880}. He then collected his results in the memoir \cite{Tissot1881},    published in 1881, in which he gives detailed proofs. In a note on p. 2 of this memoir, Tissot declares that after he published his first \emph{Comptes Rendus} notes on the subject, the statements that he gave there without proof were reproduced by A. Germain in his \emph{Trait\'e des projections des cartes g\'eographiques} \cite{Germain} and by U. Dini in his memoir \emph{Sopra alcuni punti della teoria delle superfici} \cite{Dini}. He notes that Germain and Dini gave their own proofs of these statements, which are nevertheless more complicated than those he had in mind and which he gives in the memoir \cite{Tissot1881}.
He also writes that Dini showed that the whole theory of curvature of surfaces may be deduced from the general theory that he had developed himself. In fact, Dini applied this theory to the representation of a surface on a sphere, using Gauss's methods. Tissot also says that his ideas were used in astronomy, by Herv\'e Faye, in his \emph{Cours d'astronomie de l'\'Ecole Polytechnique} \cite{Faye}.  The texts of the two \emph{Comptes rendus} notes \cite{Tissot-CR3} and \cite{Tissot-CR0} of Tissot are reproduced in the Germain's treatise \cite{Germain}.

Besides his work on geographical maps, Tissot\index{Tissot, Nicolas-Auguste (1824--1897)} wrote several papers on elementary geometry. We mention incidentally that several preeminent French mathematicians of the nineteenth and the beginning of the twentieth century published papers on this topics. We mention Serret, Catalan, Laguerre, Darboux, Hadamard and Lebesgue;  see e.g. \cite{Tissot2, Tissot3, Tissot1}.

On the title page of Tissot's\index{Tissot, Nicolas-Auguste (1824--1897)} memoir \cite{Tissot1881} (1881), the expression \emph{Examinateur \`a l'\'Ecole Polytechnique} follows his name, as he was in charge of the entrance examination. In his \emph{\'Eloge historique de Henri Poincar\'e}  \cite{Darboux-Eloge}, Darboux  relates the following episode about Tissot, examining Poincar\'e\index{Poincaré, Henri (1854--1912)}:\footnote{Poincar\'e entered the \'Ecole Polytechnique in 1873. In the French system of oral examinations, which is still in use, a student is given a question or a set of questions which he is asked to prepare while another student (who had already been given some time to prepare his questions) is explaining his solutions at the blackboard, in the same room. Thus, it is not unusual that at such an examination, some students listen to the examinations of others.}\index{Poincaré, Henri (1854--1912)}\index{Tissot indicatrix}\index{indicatrix!Tissot}
 \begin{quote}\small
 Before asking his questions to Poincar\'e, Mr. Tissot suspended the exam during 45 minutes: we thought it was the time he needed to prepare a sophisticated question. Mr. Tissot came back with a question of the Second Book of Geometry. Poincar\'e drew a formless circle, he marked the lines and the points indicated by the examiner, then, after wandering long enough in front of the blackboard, with his eyes fixed on the ground, he concluded loudly: ``It all comes down to proving the equality $AB=CD$. This is a consequence of the theory of mutual polars, applied to the two lines." Mr. Tissot interrupted him:  ``Very good, Sir, but I want a more elementary solution." Poincar\'e started wandering again, this time not in front of the blackboard, but in front of the table of the examiner, facing him, almost unconscious of his acts; then suddenly he developed a trigonometric solution. 
Mr. Tissot objected: ``I would like you to stay in Elementary Geometry." Almost immediately after that, the examiner of Elementary Geometry was given satisfaction. He warmly congratulated the examinee and announced that he deserves the highest grade.\footnote{Avant d'interroger Poincaré, M. Tissot suspendit l'examen pendant trois quarts d'heure : le temps de préparer une question raffinée, pensions-nous. M. Tissot revint avec une question du deuxième Livre de Géométrie. Poincaré dessina un cercle informe, il marqua les lignes et les points indiqués par l'examinateur; puis, après s'être promené devant le tableau les yeux fixés à terre pendant assez longtemps, conclut à haute voix: Tout revient à démontrer l'égalité $AB = CD$. Elle est la conséquence de la théorie des polaires réciproques, appliquée aux deux droites.
\\
``Fort bien, Monsieur, interrompit M. Tissot; mais je voudrais une solution plus élémentaire."  Poincaré se mit à repasser, non plus devant le tableau, mais devant la table de l'examinateur, face à lui, presque inconscient de ses actes, puis tout à coup développa une solution trigonométrique.
\\
``Je désire que vous ne sortiez pas de la Géométrie élémentaire," objecta M. Tissot, et presque aussitôt satisfaction fut donnée à l'examinateur d'élémentaires, qui félicita chaleureusement l'examiné et lui annonça qu'il avait mérité la note maxima. }
  \end{quote}

Poincaré\index{Poincaré, Henri (1854--1912)} kept a positive momory of Tissot's\index{Tissot, Nicolas-Auguste (1824--1897)} examinations. He expresses this in a letter to his mother sent on May 6, 1874, opposing them to the 10-minute examinations (known as ``colles") that he had to take regularly at the \'Ecole Polytechnique and which he said are pitiful. He writes:\footnote{\cite{Poincare-Corresp}, letter No. 62.}
 ``When I think about the exams of Tissot\index{Tissot, Nicolas-Auguste (1824--1897)} and others, I can not help but take pity of these 10 minutes little \emph{colles} where one puts in danger his future with an expression which is more or less exact or a sentence which is more or less well crafted, and where a person is judged upon infinitesimal differences."\footnote{Quand je pense aux exams de Tissot et autres, [\ldots] je ne puis m'empêcher de prendre en pitié ces petites colles de 10 minutes où on joue son avenir dans une expression plus ou moins exacte ou sur une phrase plus ou moins bien tournée et où on juge un individu sur des différences infinitésimales.}
 
 \section{On the work of Tissot on geographical maps}\label{s:work} 
 Tissot\index{Tissot, Nicolas-Auguste (1824--1897)} studied at the \'Ecole Polytechnique, an engineering school where the students had a high level of mathematical training and at a period where the applications of the techniques of differential geometry to all the domains of science were an integral part of the curriculum. 
 His work\index{Tissot, Nicolas-Auguste (1824--1897)} is part of a well-established tradition where mathematical tools are applied to the craft of map drawing. This tradition passes through the works of preeminent mathematicians such as  Ptolemy \cite{Ptolemy-geo, BJ}, Lambert \cite{Lambert-Bey, Lamb-Anmer},  Euler \cite{Euler-rep-1777, Euler-pro-1777, Euler-pro-Desli-1777} Lagrange \cite{Lagrange1779, Lagrange1779a}, Gauss \cite{Gauss-Copenhagen},  Chebyshev \cite{Cheb1, Cheb2}, Beltrami \cite{Beltrami1865}, Liouville (see the appendices to \cite{Monge}), Bonnet \cite{Bonnet-these} Darboux \cite{Darboux-Chebyshev, Darboux-Lagrange, Darboux-Tissot}, and there are others. It was known since antiquity that there exist conformal (that is, angle-preserving) projections from the sphere to the Euclidean plane.\footnote{Ptolemy, in his \emph{Geography}, works with the strereographic projection, see \cite{BJ}. See also Ptolemy's work on the \emph{Planisphere} \cite{SB}.}But it was noticed that these projections distort other quantities (length, area, etc.), and the question was to find projections that realize a compromise between these various distorsions. For instance, one question was to find the closest-to-conformal projection among the maps that are area-preserving. Hence, the idea of ``closest-to-conformal" projection came naturally. 
Among the mathematicians who worked on such problems, Tissot came closest to the notion of quasiconformality. 
 
 Let us summarize a few of his results on this subject.

 An important observation made by Tissot\index{Tissot, Nicolas-Auguste (1824--1897)} right at the beginning of his memoir \cite{Tissot1881} (p. 1) is that  
finding the most appropriate mode of projection depends on the shape of the region---and not only its size, that is, on the properties of its boundary,. Finding maps of small ``distorsion" (where, as we mentioned, this word has several possible meanings) was the aim of theoretical cartography. Tissot discovered that in order for the map to minimize an appropriately defined distortion, a certain function $\lambda$, defined by setting 
\[d\sigma^2=(1+\lambda)^2ds^2,\]
must be minimized in some appropriate sense, where $ds$ and $d\sigma$ are the line elements at the source and the target surfaces respectively. The minimality of $\lambda$ may mean, for example, that the value of the gradient of its square must be the smallest possible.

  In fact, Tissot\index{Tissot, Nicolas-Auguste (1824--1897)} studied mappings between surfaces that are more general than those between subsets of the sphere and of the Euclidean plane. He started by noting that for a given mapping between two surfaces, there is, at each point of the domain, a pair of orthogonal directions that are sent to a pair of orthogonal directions on the image surface. Unless the mapping is angle-preserving at the given point, these pairs of orthogonal directions are unique. The orthogonal directions at the various points on the two surfaces define a pair of orthogonal foliations preserved by the mapping. 
Tissot calls the tangents to these foliations \emph{principal tangents} at the given point. They correspond to the directions where the ratio of lengths of the corresponding infinitesimal line elements attains its greatest and smallest values. 

Using the foliations defined by the principal tangents, Tissot\index{Tissot, Nicolas-Auguste (1824--1897)} gave a method for finding the image of an infinitely small figure drawn in the tangent plane of the first surface. 
In particular, for a differentiable mapping, the images of infinitesimal circles are ellipses. In this case, he gave a practical way of  finding the major and minor axes of these ellipses, and he provided formulae for them. This is the theory of the 
Tissot indicatrix\index{Tissot indicatrix}\index{indicatrix!Tissot}.

  From the differential geometric point of view, the Tissot indicatrix gives information on the metric tensor obtained by pushing forward the metric of the sphere (or the spheroid) by the projection mapping.   
  
  We recall that in modern quasiconformal theory, an important parameter of a map is the quasiconformal dilatation at a point, defined as the ratio of the major axis to the minor axis of the infinitesimal ellipse which is the image of an infinitesimal circle by the map (assumed to be differentiable at the give point, so that its derivative sends circles centered at the origin in the tangent plane to ellipses).
The Tissot indicatrix\index{Tissot indicatrix}\index{indicatrix!Tissot} gives much more information than this quasiconformal dilatation, since it keeps track of (1) the direction of the great and small axes of the infinitesimal ellipse, and (2) the size of this ellipse, compared to that of the infinitesimal circle of which it is the image.

Darboux got interested in the work of Tissot\index{Tissot, Nicolas-Auguste (1824--1897)} on geography,\index{geographical map}  and in particular, in a projection described in Chapter 2 of his memoir \cite{Tissot1881}. He wrote a paper on Tissot's work \cite{Darboux-Tissot} explaining more carefully some of his results. He writes:  ``[Tissot's] exposition appeared to me a little bit confused, and it seems to me that while we can stay in the same vein, we can follow the following method [\ldots]"\footnote{Son exposition m'a paru quelque peu confuse, et il m'a semblé qu'en restant dans le même ordre d'idées on pourrait suivre avec avantage la méthode suivante.}

  Tissot\index{Tissot, Nicolas-Auguste (1824--1897)} showed then how to construct mappings that have minimal distortion. 
  
  Tissot's work was considered as very important by cartographers. The American cartographer, in his book \emph{Flattening the earth: two thousand years of map projections} \cite{Snyder-F}, published in 1997 and which is a reference in the subject, after presenting the existing  books on cartography,\index{geographical map} writes: ``Almost all of the detailed treatises presented one or two new projections, they basically discussed those existing previously, albeit with very thorough analysis. One scholar, however, proposed an analysis of distorsion that has had a major impact on the work of many twentieth-century writers on map projections. This was Tissot [\ldots]."

 Modern cartographers are still interested in the theoretical work of Tissot, see \cite{Laskowski}.

      We mentioned several  preeminent mathematicians who before Tissot worked on the theory of geographical maps.\index{geographical map}  From the more recent era, let me mention Milnor's paper titled \emph{A problem in cartography}\cite{Milnor}, published in 1969.  The reader interested in the theory of geographical maps developed by mathematicians is referred to the papers \cite{Papa-Chebyshev}, \cite{2016-Tissot} and \cite{Papa-qc} which also contain more on the work of Tissot.

 \printindex

 \end{document}